\documentclass[12pt]{article}
\usepackage{graphicx,amsmath,amsfonts,amssymb}
\usepackage[numbers]{natbib}
\usepackage[colorlinks]{hyperref}
\newcommand{\rs}[1]{{\mbox{\scriptsize \sc #1}}}

\newcommand{\sr}[1]{{\cal #1}}
 \newcommand{\dd}[1]{\mathbb{#1}}

 \newcommand{\rmn}[1]{\if#11I\else {\if#12I\hspace{-0.12ex}I\hspace{-0.85ex}\else {\if #13I\hspace{-0.16ex}I\hspace{-0.16ex}I\hspace{-1.6ex}\else I\hspace{-1.2ex}V \fi} \fi} \fi}

\newcommand{\eqn}[1]{(\ref{eqn:#1})}
\newcommand{\lem}[1]{Lemma~\ref{lem:#1}}

\newcommand{\thr}[1]{Theorem~\ref{thr:#1}}

\newcommand{\exa}[1]{Example~\ref{exa:#1}}
\newcommand{\den}[1]{Definition~\ref{den:#1}}

\newcommand{\rem}[1]{Remark~\ref{rem:#1}}

\newcommand{\sectn}[1]{Section~\ref{sect:#1}}

\newcommand{\lemt}[1]{\ref{lem:#1}}

\newcommand{\pend}{\hfill \thicklines \framebox(6.6,6.6)[l]{}}

\newenvironment{proof}{\noindent {\sc  Proof.} \rm}{\pend}
\newenvironment{proof*}[1]{\noindent {\sc  #1} \rm}{\pend}

\newtheorem{theorem}{Theorem}[section]
\newtheorem{lemma}{Lemma}[section]

\newtheorem{remark}{Remark}[section]

\newtheorem{example}{Example}[section]

\newtheorem{definition}{Definition}[section]

\newenvironment{mylist}[1]{\begin{list}{}
{\setlength{\itemindent}{#1mm}}
{\setlength{\itemsep}{0ex plus 0.2ex}}
{\setlength{\parsep}{0.5ex plus 0.2ex}}
{\setlength{\labelwidth}{10mm}}
}{\end{list}}

\begin{document}

 \title{Reversibility in Queueing Models}

\author{Masakiyo Miyazawa\\ Tokyo University of Science}
\date{April 23, 2013, R1 corrected}

\maketitle

\section{Introduction}
\label{sect:introduction}

Stochastic models for queues and their networks are usually described by stochastic processes, in which random events evolve in time. Here, time is continuous or discrete. We can view such evolution in backward for each stochastic model. A process for this backward evolution is referred to as a time reversed process. It also is called a reversed time process or simply a reversed process. We refer to the original process as a time forward process when it is required to distinguish from the time reversed process.

It is often greatly helpful to view a stochastic model from two different time directions. In particular, if some property is invariant under change of the time directions, that is, under time reversal, then we may better understand that property. {\bf A concept of reversibility} is invented for this invariance. 

A successful story about the reversibility in queueing networks and related models is presented in the book of Kelly \cite{Kell1979}. The reversibility there is used in the sense that the time reversed process is stochastically identical with the time forward process, but there are many discussions about stochastic models whose time reversal have the same modeling features. In this article, we give a unified view to them using the broader concept of the reversibility. We first consider some examples. 

\begin{mylist}{6}
\item [(i) {\bf Stationary Markov process}] (continuous or discrete time) is again a stationary Markov process under time reversal if it has a stationary distribution. This will be detailed in \sectn{Markov}. In this case, the property that transition probabilities from a given state are independent of the past history, which is called Markov property, is reversible for a stationary Markov processes.

\item [(ii) {\bf Random walk}] is a discrete time process with independent increments, and again a random walk under time reversal. That is, the time reversed process also has independent increments. Hence, the independent increment property is reversible for a random walk. A random walk is a special case of a discrete time Markov process, but does not have a stationary distribution.

\item [(iii)] {\bf Queueing model} is a system to transform an arrival stream to a departure stream. If the queueing model can be reversed in time, the roles of arrival and departure streams are exchanged, but input-output structure may be considered to be unchanged. Thus, the stream structure is reversible for the queueing model. It is notable that a stochastic process for this model does not need to have the stationary distribution. We will propose a self-reacting system for this type of models in \sectn{self-reacting}.

\item [(iv)] {\bf Birth and death process} is a continuous time Markov process taking nonnegative integers which goes up or down by 1. This process is typically used to represent population dynamics of living things. It is also used to model a queueing system, which is called a Markovian queue (see \sectn{queue}).
\end{mylist}

  These examples may be too general or too specific to get something to be useful, but suggest that there may be different levels of reversibility depending on characteristics to be of interest and on a class of stochastic processes (or models). In particular, (iv) is a special case of (i) if it has the stationary distribution.  As we will see later, this Markov chain with the stationary distribution is reversible in time, that is, the time reversed Markov chain is stochastically identical with the time forward process. Thus, (i) and (iv) are two extremes concerning the time reversed processes which have the stationary distributions.
  
  How can we define reversibility of a stochastic process for a queueing model ? Obviously, we need a time reversed dynamics. How can we construct it ? Is the stationary distribution required for this ? An answer is no because of the examples (ii) and (iii). To make the problems to be concrete, we consider discrete and continuous time Markov processes with countable state spaces, which are referred to as discrete and continuous time {\bf Markov chains}, respectively, because they are widely used in queueing and their network applications.
  
  When we are interested in the stationary distributions of those Markov chains, we have to verify or assume their existence. In this case, (i) and (iv) are extreme cases, there should be useful notion of reversibility between them. However, structural properties such as arrivals and departures should not depend on whether or not underlying Markov chains for them have the stationary distributions. We aim to device such reversibility, which will be discussed in \sectn{reversibility} under the name, {\bf reversibility in structure}. We will show how this concept arises in queueing models and their networks and how it is useful to get the stationary distributions in a unified way in the subsequent sections up to \sectn{product form}. Reversibility in structure is discussed  for other types of Markov processes in \sectn{other}, and some remarks for further study are given in \sectn{concluding}.
  
  Before we start these stories, we recall basic facts on Markov chains and their time reversals in the next section. In particular, we define a transition rate function in a slightly extended way, which will be explained.

\section{Markov chain and its time reversal}
\label{sect:Markov}

In this section, we construct time reversed processes for a discrete and continuous time Markov chains. We first consider the discrete time case because it is simpler than the continuous time case but arguments are similar. A discrete time stochastic process $\{X_{n}; n=0,1,\ldots\}$ with $S$ is called a Markov chain if, for any $n \in \dd{Z}_{+} \equiv \{0,1,2,\ldots\}$ and any $i,j,i_{0}, \ldots, i_{n-1} \in S$,
\begin{eqnarray*}
  \dd{P}(X_{n+1} = j | X_{0} =i_{0}, X_{1}=i_{1}, \ldots, X_{n-1}=i_{n-1}, X_{n}=i) = \dd{P}(X_{n+1} = j |X_{n}=i).
\end{eqnarray*}
If the right-hand side is independent of $n$, then it is said to have a time homogeneous transitions, which is denoted by $p_{ij}$. That is,
\begin{eqnarray*}
  p_{ij} = \dd{P}(X_{n+1} = j |X_{n}=i), \qquad i,j \in S.
\end{eqnarray*}
This $p_{ij}$ is called a transition probability from state $i$ to state $j$. Throughout this article, we always assume for a Markov chain to be time homogeneous. To exclude trivial exceptions, we further assume that there is a path from $i$ to $j$ for each $i,j \in S$, that is, there are states $i_{1}, i_{2}, \ldots, i_{n-1} \in S$ for some $n \ge 1$ such that
\begin{eqnarray}
\label{eqn:irreducibility 1}
  p_{ii_{1}} p_{i_{1}i_{2}} \times \cdots \times p_{i_{n-2}i_{n-1}} p_{i_{n-1}j} > 0.
\end{eqnarray}
This condition is meant that any state is commutative with any other state, which is called irreducibility.

If a probability distribution $\pi \equiv \{\pi(i)\}$ on $S$, that is, $\pi(i) \ge 0$ for all $i$ such that $\sum_{i \in S} \pi(i) = 1$, satisfies
\begin{eqnarray}
\label{eqn:stationary 1}
  \pi(j) = \sum_{i \in S} \pi(i) p_{ij}, \qquad j \in S,
\end{eqnarray}
then $\pi$ is called a stationary distribution of Markov chain $\{X_{n}\}$. Note that, $\pi(j) > 0$ for all $j \in S$ by the irreducibility. If Markov chain $\{X_{n}\}$ has a stationary distribution $\pi$ and if $X_{0}$ is subject to this $\pi$, then $X_{n}$ has the distribution $\pi$ for each $n \ge 1$. Define the time shifted process $\{X^{-\ell}_{n}; n=-\ell, -(\ell-1), \ldots\}$ by
\begin{eqnarray*}
  X^{(-\ell)}_{n} = X_{n+\ell}, \qquad n=-\ell, -(\ell-1), \ldots.
\end{eqnarray*}
Then, $\{X^{(-\ell)}_{n}; n \ge m\}$ for each $m \ge -\ell$ has the same joint distribution as $\{X_{n}; n \ge 0\}$ as long as $X_{0}$ has the stationary distribution $\pi$, and therefore, by letting $\ell \to \infty$, Markov chain $\{X^{(-\infty)}_{n}; n \in \dd{Z} \}$ is well defined, where $\dd{Z}$ is the set of all integers. Thus, we can consider the Markov chain $\{X_{n}\}$ to be a process starting from $-\infty$. This allows us to define a time reversed process $\{\tilde{X}_{n}\}$ by
\begin{eqnarray}
\label{eqn:reversed chain 1}
  \tilde{X}_{n} = X_{-n}, \qquad n \in \dd{Z}.
\end{eqnarray}

Since $\dd{P}(X_{-(n+1)}=j) = \pi(j)$, we have
\begin{eqnarray*}
  \lefteqn{\dd{P}(\tilde{X}_{0} =i_{0}, \tilde{X}_{1}=i_{1}, \ldots, \tilde{X}_{n-1}=i_{n-1}, \tilde{X}_{n}=i, \tilde{X}_{n+1} = j)}\\
  && = \dd{P}(X_{0} =i_{0}, X_{-1}=i_{1}, \ldots, X_{-(n-1)}=i_{n-1}, X_{-n}=i, X_{-(n+1)}=j)\\
  && = \dd{P}(X_{-(n+1)}=j) \dd{P}(X_{0} =i_{0}, X_{-1}=i_{1}, \ldots, X_{-(n-1)}=i_{n-1}, X_{-n}=i| X_{-(n+1)}=j)\\
  && = \pi(j) p_{ji} p_{ii_{n-1}} \times \cdots \times p_{i_{1}i_{0}},
\end{eqnarray*}
and therefore
\begin{eqnarray*}
  \lefteqn{\dd{P}(\tilde{X}_{n+1} = j|\tilde{X}_{0} =i_{0}, \tilde{X}_{1}=i_{1}, \ldots, \tilde{X}_{n-1}=i_{n-1}, \tilde{X}_{n}=i )}\\
  && = \frac { \dd{P}(\tilde{X}_{0} =i_{0}, \tilde{X}_{1}=i_{1}, \ldots, \tilde{X}_{n-1}=i_{n-1}, \tilde{X}_{n}=i, \tilde{X}_{n+1} = j)} {\dd{P}( \tilde{X}_{0} =i_{0}, \tilde{X}_{1}=i_{1}, \ldots, \tilde{X}_{n-1}=i_{n-1}, \tilde{X}_{n}=i )}\\
  && = \frac {\pi(j) p_{ji} p_{ii_{n-1}} \times \cdots \times p_{i_{1}i_{0}}} {\pi(i) p_{ii_{n-1}} \times \cdots \times p_{i_{1}i_{0}}}\\
  && = \frac {\pi(j) p_{ji}} {\pi(i)} = \dd{P}(\tilde{X}_{n+1} = j|\tilde{X}_{n}=i ).
\end{eqnarray*}
Hence, the time reversed process $\{\tilde{X}_{n}\}$ is also a Markov chain with transition probability $\tilde{p}_{ij}$ given by
\begin{eqnarray}
\label{eqn:reversed transition 1}
  \tilde{p}_{ij} = \frac {\pi(j)} {\pi(i)} p_{ji}, \qquad i,j \in S.
\end{eqnarray}

We have defined the time reversed process $\{\tilde{X}_{n}\}$ by \eqn{reversed chain 1}, but we also can define it as a Markov chain with transition probabilities $\{\tilde{p}_{ij}\}$ of \eqn{reversed transition 1}. Furthermore, for the latter definition,  it is not necessary for $\{\pi(i)\}$ to be a probability distribution, but $0 < \pi(i) < \infty$ for all $i \in S$ satisfying \eqn{stationary 1} is sufficient. This $\{\pi(i)\}$ is called a {\bf stationary measure}. This also answers why a random walk can be reversed in time because $\pi(i) = a$ for any constant $a > 0$ is its stationary measure. We refer to a measure to define $\tilde{p}$ for $p$ as in \eqn{reversed transition 1} to a reference measure (see \den{structure reversibility} and Section 3 of \cite{MiyaZwar2012} for related topics).

It is notable that we can use a reference measure other than a stationary measure. For example, we can consider a time reversed process $\{X_{n_{0}-n}; 0 \le n \le n_{0}\}$ starting at $X_{n_{0}} = i_{0}$ for fixed time $n_{0} > 0$ and state $i_{0} \in S$. In this case, the reference measure is concentrated on $i_{0}$. It can be shown that this process is again a Markov chain, but does not have time homogeneous transitions. Thus, it is different from $\{\tilde{X}_{n}\}$. As a reference measure, a stationary distribution (or measure) has the excellent feature that the transition probabilities have simple form.

Thus, the choice of a reference measure is crucial for the time reversed process. Throughout this article, we take a stationary distribution (a stationary measure for a few places) as a reference measure for a time reversed process. However, this does not mean that we also take it for reversibility in structure, which will be considered in \sectn{reversibility} because the structural property should be independent of the existence of the stationary distribution.

We next consider a continuous time Markov chain, for which irreducibility is assumed. This is quite easy because most arguments are parallel to those of a discrete time Markov chain.

A stochastic process $\{X(t)\}$ with state space $S$ is called a Markov chain if, for each $n \in \dd{Z}_{+}$, any $s>0$, any $0 < t_{1} < \ldots < t_{n} < s$ and any $i,j,i_{0}, \ldots, i_{n} \in S$,
\begin{eqnarray*}
  \lefteqn{\dd{P}(X(s+t) = j | X(0) =i_{0}, X(t_{1})=i_{1}, \ldots, X(t_{n})=i_{n}, X(s) = i)} \hspace{40ex}\\
  && = \dd{P}(X(s+t) = j |X(s)=i).
\end{eqnarray*}
We denote the right-hand side of this equation by $p_{ij}(t)$ if it is independent of $s$, which we always assume. We also assume that, for each $i,j \in S$, there exists $t > 0$ such that $p_{ij}(t) > 0$, which is called irreducibility.

Then, it is known that there exist $a(i) \ge 0$ and $q(i,j) \ge 0$ such that
\begin{eqnarray*}
  && \lim_{h \downarrow 0} \frac 1h p_{ij}(h) = q(i,j), \qquad i \ne j, i,j \in S,\\
  && \lim_{h \downarrow 0} \frac 1h (p_{ii}(h) - 1) = - a(i), \qquad i \in S.
\end{eqnarray*}
In words, $a(i)$ is a rate going out from state $i$, and $q(i,j)$ is a transition rate from $i$ to $j$ under the condition that the transition occurs. Throughout this article, we assume that $a(i)$ is finite for all $i \in S$. It is not hard to see that
\begin{eqnarray*}
  a(i) = \sum_{j \ne i, j \in S} q(i,j), \qquad i \in S,
\end{eqnarray*}
and Markov chain $\{X(t)\}$ is determined by $\{q(i,j); i,j \in S\}$. For the $X(t)$ to be well defined for all $t \ge 0$, a certain regularity condition is required (e.g., see Section 2.2 of \cite{Brem1999}). However, it is always satisfied  in queueing applications.

This continuous time Markov chain is slightly extended by adding dummy transitions by letting $q(i,i) \ge 0$. In this case, we redefine $a(i)$ as
\begin{eqnarray*}
  a(i) = \sum_{j \in S} q(i,j), \qquad i \in S.
\end{eqnarray*}
In standard text books (e.g., see \cite{Brem1999}), $q(i,i)$ is usually defined as $q(i,i) = -a(i)$. Thus, it is notable that our $q(i,i)$ is different from them. As a stochastic process, $q(i,i) > 0$ is not meaningful because it does not change the state. However, it is important for describing queueing model because it enables to rightly count arriving customers who can not enter a system. We refer to $q$ as a {\bf transition rate function}.

Since the Markov chain is irreducible, its stationary distribution $\pi$ is uniquely determined by
\begin{eqnarray}
\label{eqn:stationary 2}
  a(j) \pi(j) = \sum_{i \in S} \pi(i) q(i,j), \qquad j \in S.
\end{eqnarray}
This equation is slightly different from \eqn{stationary 1} of the discrete time Markov chain. This is just because transition rate function $q$ is not a probability distribution for given $i$. Assume that the Markov chain has the stationary distribution $\pi$, then all arguments for the discrete time Markov chain are valid for the continuous time Markov chain if we replace $p_{ij}$ by $q(i,j)$ and make small modifications like \eqn{stationary 2}. Denote the time reversed Markov chain by $\{\tilde{X}(t)\}$, then, similar to \eqn{reversed transition 1}, its transition rate function $\tilde{q}$ is given by
\begin{eqnarray}
\label{eqn:reversed q 1}
   \tilde{q}(i,j) = \frac {\pi(j)}{\pi(i)} q(j,i), \qquad i,j \in S.
\end{eqnarray}

From this equation and \eqn{stationary 2}, we immediately have the next fact.
\begin{lemma} {\rm
\label{lem:Kelly}
  A probability distribution $\pi$ on $S$ is the stationary distribution of the Markov chain with transition rate function $q$ if and only if there exists a nonnegative function $\tilde{q}$ on $S$ such that
\begin{eqnarray}
\label{eqn:Kelly 1}
 && \pi(i) \tilde{q}(i,j) = \pi(j) q(j,i), \qquad i,j \in S,\\
\label{eqn:Kelly 2}
 && \sum_{j \in S} \tilde{q}(i,j) =  \sum_{j \in S} q(i,j), \qquad i \in S.
\end{eqnarray}
}\end{lemma}

A message from this lemma is that the stationary distribution can be find if we can guess $\tilde{q}$ from the time reversed model. This idea is extensively used in the book \cite{Kell1979}, and \lem{Kelly} is called {\bf Kelly's lemma} (e.g., see \cite{ChaoMiyaPine1999}).

\section{Local balance and reversibility in structure}
\label{sect:reversibility}

From now on, we only consider the continuous time Markov chain $\{X(t)\}$ with state space $S$ and transition probabilities of $\{X(t)\}$ because the corresponding results for a discrete time Markov chain are similarly obtained. From now on, we change the notation for state from $i,j \in S$ to $x, x' \in S$ for convenience of our discussions.

We aim to define reversibility in structure for this Markov chain. For this, we first consider the strongest case that the Markov chain $\{X(t)\}$ itself is reversible. Assume that it has the stationary measure $\pi$, and define the time reversed Markov chain $\{\tilde{X}(t)\}$ with transition rate functions $\tilde{q}$. Then, $\{X(t)\}$ is stochastically identically with its time reversed process $\{\tilde{X}(t)\}$ if and only if
\begin{eqnarray*}
  q(x,x') = \tilde{q}(x',x),  \qquad x,x' \in S.
\end{eqnarray*}
From \eqn{reversed q 1}, this condition is equivalent to
\begin{eqnarray}
\label{eqn:reversibility 1}
  \pi(x) q(x,x') = \pi(x') q(x',x) , \qquad x,x' \in S.
\end{eqnarray}
This Markov chain is said to be {\bf reversible}. It is notable that \eqn{reversibility 1} yields, for fixed $x$ and the sequence of the states $x, x_{1}, \ldots, x_{n-1}, x' \in S$,
\begin{eqnarray*}
  \pi(x') = \pi(x) \frac{q(x,x_{1}) q(x_{1},x_{2}) \times \cdots \times q(x_{n-2},x_{n-1}) q(x_{n-1},x')} {q(x',x_{n-1}) q(x_{n-1},x_{n-2}) \times \cdots \times q(x_{2},x_{1}) q(x_{1},x')}.
\end{eqnarray*}
as long as the denominator in the right-hand side is positive. Hence, we have the stationary distribution $\pi$ if $\sum_{x \in S} \pi(x) < \infty$. Otherwise, this $\pi$ is a stationary measure. Thus, the reversibility enables us to get the stationary distribution. 

A typical example of a reversible Markov chain is a {\bf birth-and-death process}, which is defined in (iv) of \sectn{introduction}. Recall that its state space $S = \dd{Z}_{+}$, and the transition probability $q(x,x') > 0$ for $x,x' \in S$ if and only if $x'=x+1$, $x'=x-1 \ge 0$ or $x=x'=0$. For this Markov chain, let
\begin{eqnarray}
\label{eqn:birth-death 1}
  \pi(x) = \pi(0) \frac {q(0,1) q(1,2) \times \cdots \times q((x-1),x)} {q(x,x-1) q(x-1,x-2) \times \cdots \times q(1,0)}, \qquad x \in S,
\end{eqnarray}
then \eqn{reversibility 1} is satisfied. Hence, the birth-and-death process is reversible.

However, a Markov chain for a queueing model may not be reversible. In particular, this reversibility is too strong for queueing networks. Motivated by this, the following notion has been studied, which is weaker than reversibility.

\begin{definition}[Local balance] {\rm
\label{den:reversibility 2}
  Assume that the Markov chain $\{X(t)\}$ with transition rate function $q$ has the stationary distribution $\pi$, and let $\tilde{q}$ be the transition rate function of its time reversed process. Let ${W}$ be a set of disjoint subsets of $S \times S$. If the Markov chain satisfies
\begin{eqnarray}
\label{eqn:LB 1a}
  \sum_{(x,x') \in C} \pi(x) q(x,x') = \sum_{(x,x') \in C} \pi(x) \tilde{q}(x,x'), \qquad C \in {W},
\end{eqnarray}
or equivalently,
\begin{eqnarray}
\label{eqn:LB 1b}
  \sum_{(x,x') \in C} \pi(x) q(x,x') = \sum_{(x,x') \in C} \pi(x') q(x',x), \qquad C \in {W},
\end{eqnarray}
 then transition rate function $q$ is said to be locally balanced with respect to ${W}$.
}\end{definition}

We can weaken the local balance by using the so-called test functions. Let $\sr{M}_{+}(S \times S)$ be the set of all nonnegative valued functions on $S \times S$. For a subset $\sr{G}$ of $\sr{M}_{+}(S \times S)$, $q$ is said to be local balanced in test functions of $\sr{G}$ if
\begin{eqnarray}
\label{eqn:LB 2b}
  \sum_{x,x' \in S} \pi(x) q(x,x') f(x,x') = \sum_{x,x' \in S} \pi(x') q(x',x) f(x,x'), \qquad f \in \sr{G}.
\end{eqnarray}
Let $\sr{G}$ be the set of indicator functions $1((x,x') \in C)$ for $C \in {W}$ of \den{reversibility 2}, where $1(\cdot)$ is the function of the statement `$\cdot$' such that it equals $1$ if the statement is true and $0$ otherwise, then \eqn{LB 2b} is equivalent to \eqn{LB 1b}. Hence, \eqn{LB 2b} certainly generalizes \eqn{LB 1b}.

However, the local balance in test functions still has limitations. A problem is that transition rate functions $q$ and $\tilde{q}$ are directly related through balance equations. This limits flexibility of $\tilde{q}$ to be chosen. Furthermore, it also requires that $q$ has the stationary distribution. Instead of balance equations and stationary distribution, we use a set of pairs of transition rate functions $q$ and measure $\pi$ on $S$ to support $q$ for a weaker sense of reversibility, where $\pi$ is said to be a measure to support $q$ or a {\bf supporting measure} of $q$ if $\pi(x) \ge 0$ for all $x \in S$, and, for each  $x\in S$, $\pi(x) > 0$ if $q(x',x) > 0$ for some $x' \in S$.

Let $\sr{T}$ be the set of all pairs $(q,\pi)$ of transition rate function $q$ on $S \times S$ and measure $\pi$ on $S$ to support $q$. Note that $\sum_{x' \in S} q(x,x') < \infty$ for all $x \in S$ is always assumed, but $\pi$ is not necessarily the stationary distribution of $q$. Furthermore, a Markov chain with transition rate function $q$ may not be irreducible. We then define the following reversibility.

\begin{definition}[Reversible in structure] {\rm
\label{den:structure reversibility}
  Let $\sr{Q}$ be a subset of $\sr{T}$. $\sr{Q}$ is said to be reversible in structure if $(q,\pi) \in \sr{Q}$ implies $(\tilde{q},\pi) \in \sr{Q}$ for all $(q,\pi) \in \sr{Q}$, where
\begin{eqnarray*}
  \tilde{q}(x,x') = \left\{\begin{array}{ll}
  \frac {\pi(x')} {\pi(x)} q(x',x), \quad & \pi(x) >0, x,x' \in S,\\
  0 & \pi(x) = 0, x,x' \in S.
  \end{array} \right.
\end{eqnarray*}
 This $\pi$ is called a {\bf reference measure} to get $\tilde{q}$ from $q$. In particular, a Markov chain with transition rate function $q$ is said to be reversible in structure with respect to $\sr{Q}$ if $(q,\pi) \in \sr{Q}$ implies $(\tilde{q}, \pi) \in \sr{Q}$ for all $(q,\pi) \in \sr{Q}$.
}\end{definition}

\begin{remark} {\rm
\label{rem:reversible in structure 1}
  If $\pi$ is the stationary distribution of $q$, then $(q, \pi) \in \sr{T}$, and $\pi$ is also the stationary distribution of $\tilde{q}$. However, the reversibility in structure does not require for $\pi$ to be the stationary distribution or stationary measure of $q$ for $(q, \pi) \in \sr{Q}$. In fact, the notion of reversibility in structure is independent of their existence.
}\end{remark}

Note that $\tilde{\tilde{q}}$ may not be identical with $\tilde{q}$. For example, this is the case when $q(x',x) > 0$ for $\pi(x) > 0$ and $\pi(x') = 0$. However, $\tilde{\tilde{\tilde{q}}}$ is $\tilde{q}$, so it is easy to see that $(q, \pi) \in \sr{Q}$ is reversible in structure concerning $\sr{Q}$ if and only if $(\tilde{q}, \pi) \in \sr{Q}$ does so. This reversibility obviously includes the local balance in test functions with respect to each $\sr{G}$. For this, it is sufficient to let $\sr{Q}$ be the set of $(q,\pi)$ such that \eqn{LB 2b} holds for all $g \in \sr{G}$ and $q$ has the stationary distribution $\pi$. However, $\sr{Q}$ of \den{structure reversibility} may be too flexible for making $q$ to be specific.

Because we are interested in queueing models, we have arrivals and departures to a system, and their roles are exchanged under time reversal. This suggests to partially decompose $q$ into a family of transition functions $q_{u}$ for $(q_{u}, \pi) \in \sr{Q}$ with $u \in U$, where $U$ is a finite set, in such a way that
\begin{eqnarray}
\label{eqn:bounded 1}
  \sum_{u \in U} q_{u}(x,x') \le q(x,x'), \qquad x,x' \in S,
\end{eqnarray}
 and to define a one-to-one mapping $\Gamma$ from $U$ to itself, that is, $\Gamma$ is a permutation on $U$. We refer to $\{q_{u}; u \in U\}$ satisfying \eqn{bounded 1} as a {\bf sub-transition family} of $q$.

For example, if $u$ represents arrivals, that is, $q_{u}$ is the transition function for arrivals, then $\Gamma(u)$ may be considered to represent departures. Since arrivals are changed to departures under time reversal, we define the transition function $\tilde{q}_{\Gamma(u)}$ for $u \in U$ using supporting measure $\pi$ of $q_{u}$ by
\begin{eqnarray}
\label{eqn:dual q u}
  \tilde{q}_{\Gamma(u)}(x,x') = \left\{\begin{array}{ll}
  \frac {\pi(x')} {\pi(x)} q_{u}(x',x), \quad & \pi(x) >0, x,x' \in S,\\
  0 & \pi(x) = 0, x,x' \in S,
  \end{array} \right.
\end{eqnarray}
that is, $\tilde{q}_{\Gamma(u)}$ is obtained from $q_{u}$ by the reference measure $\pi$. We may consider $\tilde{q}_{\Gamma(u)}$ to be the transition rate function for departures under time reversal with respect to $\pi$. Note that \eqn{dual q u} can be used to define $\tilde{q}_{u}$ if we use $\Gamma^{-1}(u)$ instead of $u$. These facts motivate us to define the following reversibility in structure.

\begin{definition}[$\Gamma$-reversible in structure] {\rm
\label{den:UG-structure reversibility}
  For a finite set $U$ and a permutation $\Gamma$ on $U$, let $\sr{Q}_{u}$ be a subset of $\sr{T}$ for each for $u \in U$, then $\{\sr{Q}_{u}; u \in U\}$ is said to be $\Gamma$-reversible in structure if $(q_{u},\pi) \in \sr{Q}_{u}$ implies $(\tilde{q}_{u},\pi) \in \sr{Q}_{u}$ for all $u \in U$, where
\begin{eqnarray*}
  \tilde{q}_{u}(x,x') = \left\{\begin{array}{ll}
  \frac {\pi(x')} {\pi(x)} q_{\Gamma^{-1}(u)}(x',x), \quad & \pi(x) >0, x,x' \in S.\\
  0 & \pi(x) = 0, x,x' \in S,
  \end{array} \right.
\end{eqnarray*}
In particular, a Markov chain which has transition rate function $q$ and sub-transition family $\{q_{u}; u \in U\}$ is said to be $\Gamma$-reversible in structure with respect to $\{\sr{Q}_{u}; u \in U\}$ if $(q_{u},\pi) \in \sr{Q}_{u}$ implies $(\tilde{q}_{u},\pi) \in \sr{Q}_{u}$ for all $u \in U$ and all supporting measures $\pi$ on $S$.
}\end{definition}

\begin{remark} {\rm
\label{rem:UG-structure reversibility}
  The $\Gamma$-reversibility in structure of $\{\sr{Q}_{u}; u \in U\}$ is stronger notion than a Markov chain with $q$ and $\{q_{u}; u \in U\}$ to be $\Gamma$-reversible in structure with respect to $\{\sr{Q}_{u}; u \in U\}$. For the latter, $\{\sr{Q}_{u}; u \in U\}$ itself is not necessarily $\Gamma$-reversible in structure. However, as we will see \lem{reversibility in structure} below, the $\Gamma$-reversibility in structure of $\{\sr{Q}_{u}; u \in U\}$ is convenient to construct a class of models satisfying such reversibility.
}\end{remark}

We simplify the conditions for $\Gamma$-reversible in structure. For this, we partition the set $U$ into permutation invariant sets. For each $u \in U$, let $U_{u}$ be the subset of $U$ such that $u \in U_{u}$ and $v \in U_{u}$ implies $\Gamma(v) \in U_{u}$. Then, either $U_{u} = U_{v}$ or $U_{u} \cap U_{v} = \emptyset$ holds for $u, v \in U$, and therefore there is the subset $U'$ such that $U$ is partitioned into disjoint sets $U_{u}$ for $u \in U'$. These sets $U_{u}$ are said to be {\bf irreducible under permutation $\Gamma$}. Obviously, for $\Gamma$-reversibility in structure, it is sufficient to check it for each irreducible set. The following lemma simplifies the conditions of $\Gamma$-reversibility in structure for an irreducible set.
\begin{lemma} {\rm
\label{lem:reversibility in structure}
  Let $\Gamma$ be a permutation on a finite set $U$. Then, $\{\sr{Q}_{u}; u \in U\}$ is $\Gamma$-reversible in structure if and only if 
\begin{eqnarray}
\label{eqn:G reversible condition 1}
  \sr{Q}_{\Gamma(u)} = \left\{(\tilde{q}_{\Gamma(u)},\pi) \in \sr{T}; (q_{u}, \pi) \in \sr{Q}_{u} \right\}, \qquad u \in U.
\end{eqnarray}
In particular, if $U$ is irreducible under permutation $\Gamma$ and if $\sr{Q}_{u}$ is sequentially defined by \eqn{G reversible condition 1} for all $u \in U$ from $u_{0} \in U$ in such a way that $u_{i} = \Gamma(u_{i-1})$ for $i=1,2, \ldots, m$ and $u_{0} =u_{m}$ for some positive integer $m$, then $\sr{Q}_{0} = \sr{Q}_{m}$ implies that $\{\sr{Q}_{u}; u \in U\}$ is $\Gamma$-reversible in structure.

\begin{proof}
Since \eqn{G reversible condition 1} is equivalent to that $(q_{u},\pi) \in \sr{Q}_{u}$ for all $u \in U$ if only if $(\tilde{q}_{\Gamma(u)},\pi) \in \sr{Q}_{\Gamma(u)}$ for all $u \in U$, the necessity and sufficiency of \eqn{G reversible condition 1} are immediate from the definition of $\Gamma$-reversibility in structure for $\{\sr{Q}_{u}; u \in U\}$. The last part of the lemma is a direct consequence of the condition that $U$ is irreducible under permutation $\Gamma$.
\end{proof}
}\end{lemma}

We will show how $\Gamma$-reversibility in structure and \lem{reversibility in structure} work for queueing and their network models in the subsequent sections. In the rest of this section, we briefly discuss a Poisson process, which is frequently used in queueing applications as a counting process of customers. We are interested in how it arises in a Markov chain describing a queueing model and how it is connected to structure reversibility.
  
\begin{definition}[Poisson process] {\rm
\label{den:Poisson 1}
  A counting process $N(t)$ is called a Poisson process with rate $\lambda>0$ if it satisfies the following conditions.
\begin{itemize}
\item [(a1)] It has independent increments, that is, for any $n \ge 2$ and any time sequence $t_{0} \equiv 0< t_{1} < t_{2} < \ldots < t_{n}$, the increments $N(t_{\ell}) - N(t_{\ell-1})$, $\ell=1,2,\ldots,n$, are independent.
\item [(a2)] It has stationary increments, that is, for any $s,t>0$, the distribution of $N(s+t) - N(s)$ is independent of $s$.
\item [(a3)] It is simple, that is, the increments $N(t) - N(t-)$ at time $t$ is at most $1$.
\item [(a4)] $\lambda \equiv \dd{E}(N(1)) < \infty$.
\end{itemize}  
}\end{definition}

  It is not hard to see that this definition is equivalent to that the inter-occurrence times for counting are independent and identically distributed subject to the exponential distribution with mean $1/\lambda$. 
  
  We here note two important facts related to reversibility. Let $q$ be the transition rate function of the continuous time Markov chain with the stationary distribution $\pi$ such that $\pi(x) > 0$ for all $x \in S$. Let $q_{*}$ be a transition rate function satisfying
\begin{eqnarray}
\label{eqn:bounded 2}
  q_{*}(x,x') \le q(x,x'), \qquad x,x' \in S.
\end{eqnarray}
Let $N_{*}(t)$ be the number of transitions of this Markov chain in the time interval $(0,t]$ generated by $q_{*}$. $N_{*}$ is called a counting process of $q_{*}$.
\begin{lemma} {\rm
\label{lem:Poisson 1}
The counting process $N_{*}$ is the Poisson process with rate $\lambda > 0$ whose occurrences are independent of the state of Markov chain just before their counting instants if and only if
\begin{eqnarray}
\label{eqn:Poisson 1}
  \sum_{x' \in S} q_{*}(x,x') = \lambda, \qquad x \in S.
\end{eqnarray}
}\end{lemma}

\lem{Poisson 1} easily follows from the following facts. The sojourn time at a given state of the Markov chain has an exponential distribution, the inter-occurrence times of $N_{*}$ are exponentially distributed, and mutually independent by the Markov property (see \cite{ChaoMiyaPine1999} for a complete proof).

For $q_{*}$ and $\pi$, recalling that $\pi(x) > 0$ for all $x \in S$, define $\tilde{q}_{*}$ as
\begin{eqnarray*}
  \tilde{q}_{*}(x,x') = \frac {\pi(x')} {\pi(x)} q_{*}(x',x), \quad & x,x' \in S.
\end{eqnarray*}
and let $\tilde{N}_{*}$ be the point process generated by $\tilde{q}_{*}$. A question is when $\tilde{N}_{*}$ is a Poisson process. This is equivalent for $\tilde{q}_{*}$ to satisfy \eqn{Poisson 1} instead of $q_{*}$. Namely,
\begin{eqnarray}
\label{eqn:Poisson 2}
  \sum_{x' \in S} \pi(x') q_{*}(x',x) = \pi(x) \lambda, \qquad x \in S.
\end{eqnarray}
This is intuitively clear because $\tilde{N}_{*}$ is the time reversed process of $N_{*}$ under $\pi$, and a Poisson process is unchanged by time reversal.

We also can answer this question in terms of reversibility in structure. For this, let $\sr{Q}_{+}$ and $\sr{Q}_{-}$ be the sets of  $\sr{T}$ satisfying
\begin{eqnarray*}
  \sr{Q}_{-} = \{(\tilde{g}_{*}, \sigma) \in \sr{T}; (g_{*}, \sigma) \in \sr{Q}_{+} \},
\end{eqnarray*}
$\tilde{g}(x,x') = \frac {\sigma(x')} {\sigma(x)} g_{*}(x',x)$ for $x,x' \in S$ satisfying $\sigma(x) > 0$ and $\tilde{g}_{*}(x,x') = 0$ otherwise. Then, letting $U = \{+.-\}$, $g_{+} = g_{*}$, $g_{-} = \tilde{g}_{*}$ and $\Gamma$ is the permutation on $U$, we can see that $\{\sr{Q}_{+}, \sr{Q}_{-}\}$ is reversible in structure by the last half of \lem{reversibility in structure}. In particular, $(q_{*},\pi) \in Q_{+}$ if and only if $(\tilde{q}_{*}, \pi) \in Q_{-}$. Thus, we have the following answer from \lem{Poisson 1}. 

\begin{lemma} {\rm
\label{lem:Poisson 2}
Assume that the Markov chain $\{X(t)\}$ has the stationary distribution $\pi$. The counting process $\tilde{N}_{*}$ is the Poisson process with rate $\lambda > 0$ whose occurrences are independent of the state of Markov chain just after their counting instants if and only if either $(q_{*}, \pi) \in \sr{Q}_{-}$, or $q$ with sub-transitions $\{q_{*}, \tilde{q}_{*}\}$ is $\Gamma$ reversible in structure with respect to $\{\sr{Q}_{+}, \sr{Q}_{-}\}$.
}\end{lemma}

\section{Queue and reacting system}
\label{sect:queue}

We now describe queueing models by a continuous time Markov chain. As a toy example, we first consider a single server queue. Assume that customers arrive according to the Poisson process with rate $\lambda > 0$, customers are served in the manner of first-come and first-served with independently and exponentially distributed service times with mean $1/\mu$ ($\mu > 0$), which are also independent of the arrival process. This model is called an {\bf $M/M/1$ queue}.

Let $X(t)$ be the number of customers at time $t$ of the $M/M/1$ queue. Since the remaining times to the next arrivals and to the service completions of a customer being served are independent and subject to exponential distributions, $X(t)$ is a continuous time Markov chain with state space $S = \{0,1,\ldots\}$, and its transition rate function $q$ is given by
\begin{eqnarray*}
  q(x,x') = \left\{\begin{array}{ll}
  \lambda, \quad & x'=x+1,\\
  \mu, \quad & x'=x-1\ge 0,\\
  0, & \mbox{otherwise},
  \end{array} \right. \qquad x,x' \in S.
\end{eqnarray*}
  Hence, $\{X(t)\}$ is a special case of a birth-and-death process, and therefore reversible if it has the stationary distribution. From \eqn{birth-death 1}, the stationary measure $\pi$ is given by
\begin{eqnarray*}
  \pi(x) = \pi(0) \rho^{x}, \qquad x \in S,
\end{eqnarray*}
where $\rho = \lambda/\mu$, and satisfies \eqn{reversibility 1}. Hence, the stationary distribution exists if and only if $\rho < 1$.

The $M/M/1$ queue is a simple model, but it has key ingredients of reversibility. For this, we introduce transition functions for arrivals and departures. Let
\begin{eqnarray*}
  q_{a}(x,x') = \lambda 1(x'=x+1), \qquad q_{d}(x,x') = \lambda 1(x'=x-1 \ge 0).
\end{eqnarray*}
Then,
\begin{eqnarray}
\label{eqn:Poisson arrival 1}
 && \sum_{x' \in S} q_{a}(x,x') = \lambda, \qquad x \in S,\\
\label{eqn:Poisson departure 1}
 && \sum_{x' \in S} \pi(x') q_{d}(x',x) = \pi(x) \lambda, \qquad x \in S.
\end{eqnarray}
Hence, by Lemmas \lemt{Poisson 1} and \lemt{Poisson 2}, both $N_{a}$ and $N_{d}$ are the Poisson processes with the same rate $\lambda$. Also, by these lemmas, $N_{a}$ is independent of the past history of the Markov chain while $N_{d}$ is independent of its future history. These facts are also valid for the $M/M/s$ queue which increases the number of servers in the $M/M/1$ queue from one to $s$ since the Markov chain for the queue is also a birth-and-death process. These facts are known as {\bf Burke's theorem} \cite{Burk1968}. 
  
  We now closely look at the $M/M/1$ queue and Burke's theorem, then we can see that Poisson arrivals and Poisson departures occur as long as \eqn{Poisson arrival 1} and \eqn{Poisson departure 1} are satisfied, and it may not be necessary for the Markov chain $\{X(t)\}$ to be a birth-and-death process. In this situation, we need to have a larger state space $S$ and to suitably choose $(q_{a},\pi) \in \sr{T}$ for arrival and $(q_{d},\pi) \in \sr{T}$ for departures, where $\pi$ is the stationary distribution of $\{X(t)\}$. Another issue is that we may not be able to define the arrival process a priori because it may be a superposition of departure processes of other queues in network applications.
  
  For resolving these issues, we separate the exogenous arrival process from the queueing system, and introduce the following formulation according to \cite{ChaoMiyaPine1999}. Let $S$ be the state space of the queueing system removing the arrival process. Assume the following dynamics.
\begin{itemize}
\item [(b1)] A customer under state $x$ departs with rate $q_{d}(x, x')$ changing $x$ to $x'$.
\item [(b2)] An arriving customer who finds state $x \in S$ changes it to $x' \in S$ with probability $p_{a}(x, x')$, where $\sum_{x' \in S} p_{a}(x, x') = 1$ for each $x \in S$.
\item [(b3)] The system state $x$ changes to $x'$ with rate $q_{\rs{i}}(x, x')$ without departure. This state transition is said to be internal.
\end{itemize}

This model is fairly general as a queueing system, but nothing to do with a customer arrival process. We refer to this model as a {\bf reacting system}.
  
  We fed Poisson arrivals with rate $\alpha > 0$ to this reacting system, and define transition rate functions by
\begin{eqnarray}
\label{eqn:reacting system 1}
  q(x,x') = \alpha p_{a}(x, x')  + q_{d}(x, x')+ q_{\rs{i}}(x, x'), \qquad x, x' \in S.
\end{eqnarray}
  Denote the Markov chain with transition rate function $q$ by $\{X(t)\}$. This Markov chain represents the dynamics given by (b1), (b2) and (b3), where all actions except for (b1) occur with exponentially distributed inter-occurrence times.

Assume that this Markov chain has the stationary distribution $\pi$. Then, the transitions rate function $\tilde{q}$ of the time reversed Markov chain $\{\tilde{X}(t)\}$ is given by \eqn{reversed q 1}. A question is under what conditions $\{\tilde{X}(t)\}$ represents a reacting system with Poisson arrivals.

This problem can be stated using $\Gamma$-reversibility in structure. For this, define $q_{a}$ as
\begin{eqnarray*}
  q_{a}(x,x') = \alpha p_{a}(x, x'), \qquad x, x' \in S,
\end{eqnarray*}
and let $U = \{a,d\}$ and let $\Gamma$ be the permutation on $U$ such that $\Gamma(a) = d$. Define $\sr{Q}_{a}$ as
\begin{eqnarray}
\label{eqn:Q a}
  \sr{Q}_{a} = \Big\{(g_{a},\sigma) \in \sr{T}; \exists c > 0, \sum_{x' \in S} g_{a}(x,x') = c \Big\}.
\end{eqnarray}
Then, $(q_{a},\pi) \in \sr{Q}_{a}$. In other words, $\{X(t)\}$ is a reacting system if and only if $(q_{a}, \pi) \in \sr{Q}_{a}$. The question is what is $\sr{Q}_{d}$ for $q $ with $\{q_{a}, q_{d}\}$ to be $\Gamma$-reversible in structure with respect to $\{\sr{Q}_{a}, \sr{Q}_{d}\}$. The answer is immediate from \lem{reversibility in structure} because $U$ is irreducible. Namely, we must have
\begin{eqnarray*}
  \sr{Q}_{d} = \{(\tilde{g}_{d}, \sigma) \in \sr{T}; (g_{a}, \sigma) \in \sr{Q}_{a} \},
\end{eqnarray*}
where $\tilde{g}_{d}(x,x') = \frac {\sigma(x')} {\sigma(x)} g_{a}(x',x)$ if $\sigma(x) > 0$ and it vanishes otherwise by the defineition $\Gamma$. This $\sr{Q}_{d}$ can be written as $\{(g_{d}, \sigma) \in \sr{T}; (\tilde{g}_{a}, \sigma) \in \sr{Q}_{a} \}$, and therefore
\begin{eqnarray}
\label{eqn:Q d}
  \sr{Q}_{d} = \Big\{(g_{d}, \sigma) \in \sr{T}; \exists c > 0, \sum_{x' \in S} \sigma(x') g_{d}(x',x) = c \, \sigma(x)\Big\},
\end{eqnarray}
since $\tilde{g}_{a}(x,x') = \frac {\sigma(x')} {\sigma(x)} g_{d}(x',x)$ for $\sigma(x) > 0$.

Let $q_{d}$ be a sub-transition rate function of $q$, and assume that $q$ has its stationary distribution $\pi$, then the condition for $(q_{d}, \pi) \in \sr{Q}_{d}$ of \eqn{Q d} can be written as, for some $\beta > 0$,
\begin{eqnarray}
\label{eqn:Poisson departure 2}
  \sum_{x' \in S} \pi(x') q_{d}(x', x) = \beta \pi(x) , \qquad x \in S,
\end{eqnarray}
which is called {\bf quasi-reversibility} in the literature (see, e.g., \cite{ChaoMiyaPine1999}). Thus, we have the following theorem with help of Lemmas \lemt{reversibility in structure} and \lemt{Poisson 2}.

\begin{theorem} {\rm
\label{thr:reversibility 1}
  Define $\sr{Q}_{a}$ and $\sr{Q}_{d}$ as \eqn{Q a} and \eqn{Q d}, respectively, then $\{\sr{Q}_{a}, \sr{Q}_{d}\}$ is $\Gamma$-reversible in structure. Let $\{X(t)\}$ be the Markov chain with transition rate functions \eqn{reacting system 1} for the reacting system with Poisson arrivals, and assume that this Markov chain has the stationary distribution $\pi$. Then, the following conditions are equivalent.
\begin{itemize}
\item [(c1)] $\{\tilde{X}(t)\}$ is also a reacting system. 
\item [(c2)] $\{X(t)\}$ with $\{q_{a}, q_{b}\}$ is $\Gamma$-reversible in structure with respect to $\{\sr{Q}_{a}, \sr{Q}_{d}\}$. 
\item [(c3)] $\{X(t)\}$ is quasi-reversible, that is, \eqn{Poisson departure 2} holds for some $\beta > 0$.
\item [(c4)] The departure process $N_{d}$ from the reacting system is the Poisson process with rate $\beta$ which is independent of the future evolution of the Markov chain $\{X(t)\}$.
\end{itemize}
}\end{theorem}

\begin{example}[$M/M/1$ batch service queue] {\rm
\label{exa:reacting system}
Modify the $M/M/1$ queue in such a way that a batch of customers depart at the service completion instant subject to the following rule. All customers depart if a requested batch size is greater than the number of customers in system at the moment. Otherwise, the requested number of customers depart. It is assumed that the batch sizes to be requested are independently and identically distributed with a finite mean. We call this model $M/M/1$ batch service queue. This queue does not have Poisson departures if all batch departures are counted as departures, but it does if we only count full batches, that is, the departing batches whose sizes meet those to be requested. Thus, this reacting system is reversible in structure if full batches are only counted as departures. Details on this model can be found in Section 2.6 of \cite{ChaoMiyaPine1999}.
\pend  
}\end{example}

In \thr{reversibility 1}, we have assumed the Poisson arrivals. It may be questioned what will occur if the arrival process is not Poisson. To answer it, we consider an external source to produce arrivals as another reacting system that is independent of the original reacting system, where departures from the latter system become arrivals to the external source but do not cause any state changes of the source. Thus, we have the closed system in which the two reacting systems are  cyclically connected. If the time reversed model has the same structure, it is not hard to see that the external source produces Poisson arrivals and departures from the original system are also subject to Poisson. Hence, the Poisson arrivals arise under the reversibility in structure of this cyclic model.

A reacting system with Poisson arrivals can model not only a single node queue but also a network of queues if the state space $S$ is appropriately chosen. A problem is whether the reversibility enables us to get the stationary distribution. Obviously, this is easier for a smaller system. This motivates us to consider reacting systems as nodes of a network queue and to construct the stationary distribution of the network from those of reacting systems of Poisson arrivals. This program will be executed in \sectn{product form}.

Another question on the reacting system is whether or not we can consider arrival processes depending on the state of the reacting system. We solve this problem incorporating exogenous arrivals as one of departures suitably extending the state space $S$. This will be done in \sectn{self-reacting}.

Before looking at those problems, we will discuss a concrete example for the quasi-reversibility.

\section{Symmetric service}
\label{sect:symmetric}

The examples in \sectn{queue} are reversible in structure, but have no internal state transition. In this section, we introduce the reacting system which have internal state transitions. For this, we classify customers by types, and split the transition functions by them.

Let $T$ be the set of customer types, and classify arrivals and departures by customer types. We split $p_{a}$, $q_{d}$, $\alpha$ and $\beta$ into $p_{au}$, $q_{du}$, $\alpha_{u}$ and $\beta_{u}$ with $u \in T$, respectively. Thus, the transition rate of the Markov chain is changed to
\begin{eqnarray*}
  q(x,x') = \sum_{u \in T} \Big(\alpha_{u} p_{au}(x, x') + q_{du}(x, x') \Big) + q_{\rs{i}}(x, x'), \qquad x, x' \in S.
\end{eqnarray*}
Similar changed are made for the time reversed model, and all the arguments are parallel for the reversibility. For example, the quasi-reversibility condition \eqn{Poisson departure 2} is changed to
\begin{eqnarray}
\label{eqn:Poisson departure 3}
  \sum_{x' \in S} \pi(x') q_{du}(x', x) = \beta_{u} \pi(x) , \qquad x \in S, u \in T.
\end{eqnarray}
Let $U = \{au, du; u \in T\}$, and let $\Gamma$ be the permutation such that $\Gamma(au) = du$ for all $u \in T$. We define $\sr{Q}_{au}$ and $\sr{Q}_{du}$ similarly to $\sr{Q}_{a}$ of \eqn{Q a} and $\sr{Q}_{d}$ of \eqn{Q d}, respectively.

To describe service discipline, we assume that there are infinitely many service positions numbered by $1,2,\ldots$, and customers in system are allocated to positions $1$ to $n$ when $n$ customers are in system. We further assume that each $x$ in the state space $S$ has the following form when $n$ customers are in system.
\begin{eqnarray*}
  x = (n, \{(u_{\ell}, w_{\ell}); \ell=1,2,\ldots,n\}), \qquad u_{\ell} \in T, w_{\ell} \in \{1,2,\ldots\},
\end{eqnarray*}
where $u_{\ell}$ and $w_{\ell}$ are the type and remaining service stage, respectively, of a customer in position $\ell$. Let $k_{u}$ be the service stage of a newly arriving customer of type $u$, and let $\eta_{uj}$ be the completion rate of the stage $j$, that is, the sojourn time at the stage $j$ is exponentially distributed with mean $1/\eta_{uj}$. We assume the following dynamics under the state $x$.
\begin{itemize}
\item [(d1)] A type $u$ customer arrives and chooses position $i$ among $\{1,2,,\ldots,n+1\}$ with probability $\delta(i,n+1)$, and customers in positions $i, i+1,\ldots, n$ move to $i+1, i+2,\ldots, n+1$. That is, $x$ is changed to $x'$:
\begin{eqnarray*}
  x' = (n+1, \{(u_{\ell}, w_{\ell}); \ell=1,2,\ldots,i\}, (u,k_{u}), \{(u_{\ell}, w_{\ell}); \ell=i+1,i+2,\ldots,n\}).
\end{eqnarray*}
\item [(d2)] A type $u_{i}$ customer in position $i$ with $w_{i} = 1$ departs with rate $\gamma(i,n) \eta_{u_{i}1}$, and customers in positions $i+1, i+2,\ldots, n$ move to $i, i+1,\ldots, n-1$. That is, $x$ is changed to $x'$:
\begin{eqnarray*}
  x' = (n-1, \{(u_{\ell}, w_{\ell}); \ell=1,2,\ldots,i-1\}, \{(u_{\ell}, w_{\ell}); \ell=i+1,i+2,\ldots,n\}).
\end{eqnarray*}
\item [(d3)] A type $u_{i}$ customer in position $i$ with $w_{i} \ge 2$ decreases its remaining service stage by 1 with rate $\gamma(i,n) \eta_{u_{i} w_{i}}$, and customers in positions $\ell \ne i$ are unchanged. That is, $x$ is changed to $x'$:
\begin{eqnarray*}
  x' = (n, \{(u_{\ell}, w_{\ell}); \ell=1,2,\ldots,i-1\}, (u_{i}, w_{i}-1), \{(u_{\ell}, w_{\ell}); \ell=i+1,i+2,\ldots,n\}).
\end{eqnarray*}
\end{itemize}

Note that (d1), (d2) and (d3) uniquely specify $p_{a}$, $q_{d}$ and $q_{\rs{i}}$. This reacting system is fairly general, but hard to get the stationary distribution even for Poisson arrivals. Thus, we consider the simpler situation assuming that, for each $n \ge 1$,
\begin{eqnarray}
\label{eqn:symmetric service 1}
  \phi(n) \delta(\ell,n) = \gamma(\ell,n), \qquad \ell=1,2,\ldots n,
\end{eqnarray}
where $\phi(n) = \sum_{\ell=1}^{n} \gamma(\ell,n)$. Service discipline satisfying this condition is referred to as {\bf symmetric service} (see, e.g., \cite{Kell1979}). This service discipline includes processor sharing and preemptive last-come last-served, but can not be first-come and first-served except for exponentially distributed service times.

It is well known that the reacting system with Poisson arrivals and symmetric service, which is called a symmetric queue, has the stationary distribution $\pi$:
\begin{eqnarray*}
  \pi(n, \{(u_{\ell}, w_{\ell}); \ell=1,2,\ldots,n\}) = c^{-1} \prod_{\ell=1}^{n} \frac {\alpha_{u_{\ell}}}{\phi(\ell) \eta_{u_{\ell}w_{\ell}}}, \qquad 1 \le w_{\ell} \le k_{u_{\ell}},
\end{eqnarray*}
 if it exists, where $c$ is the normalizing constant (e.g., see \cite{ChaoMiyaPine1999}). Furthermore, the quasi-reversibility \eqn{Poisson departure 3} holds, and therefore $(\tilde{q}_{au},\pi) \in \sr{Q}_{au}$. Similarly, $(\tilde{q}_{ud}, \pi) \in \sr{Q}_{du}$, where 
\begin{eqnarray*}
  \tilde{q}_{du}(x,x') = \frac {\pi(x')} {\pi(x)} \alpha p_{au}(x', x), \qquad x, x' \in S.
\end{eqnarray*}
Thus, the symmetric queue is $\Gamma$-reversible in structure with respect to $\{\sr{Q}_{au}, \sr{Q}_{du}; u \in T\}$. Furthermore, we can verify that
\begin{eqnarray*}
  q_{\rs{i}}(x,x') = \tilde{q}_{\rs{i}}(x,x'), \qquad x,x' \in S,
\end{eqnarray*}
These facts are particularly useful in incorporating  symmetric queues as nodes of a network.

In the symmetric service, the service amount of type $u$ customer is the sum of $k_{u}$ independent exponentially distributed random variables with means $1/\eta_{uk_{u}}, 1/\eta_{u(k_{u}-1)}, \ldots, 1/\eta_{u1}$. Hence, if we let $\eta_{uj} = k_{u} \zeta_{u}$, then the service amount distribution of type $u$ customer is $k_{u}$-order Erlang distribution with mean $\zeta_{u}$. Furthermore, if type $u$ is split to finitely many sub-types and if each sub-type has finitely many stages, then the amount of work for type $u$ customer is subject to the mixture of Erlang distributions. Thus, we can approximate any service amount distribution by this finitely many stage model.

\section{Self-reacting system and balanced departure}
\label{sect:self-reacting}

We modify the reacting system so that it includes the exogenous source as a part of the system, and departures are immediately transferred to arrivals. We allow to have different types for arrivals and departures. Let $T$ be the set of types. In this section, an element of $T$ may represent characteristics other than a type of customer. For example, it may represent batch sizes of arrivals and departures and their vectors. We slightly change the dynamics (b1) and (b2) while keeping (b3).
\begin{itemize}
\item [(e1)] An entity $u$ is released under state $x$ with rate $q_{du}(x, y)$ changing state $x$ to $y$.
\item [(e2)] A released entity $u$ is transferred to entity $u'$ under state $y$ with probability $r(u,y,u')$, where
\begin{eqnarray*}
  \sum_{u \in T} r(u,y,u') = 1, \qquad u \in T, y \in S.
\end{eqnarray*}
\item [(e3)] A transferred entity $u'$ under state $y \in S$ changes the state from $y$ to $x' \in S$ with probability $p_{au'}(y, x')$, where $\sum_{x' \in S} p_{au'}(y, x') = 1$ for each $u' \in T$ and $y \in S$.
\end{itemize}

We assume that (e2) and (e3) instantaneously follows after (e1). Then, we can define the transition rate function $q$ of the Markov chain $\{X(t)\}$ as
\begin{eqnarray}
\label{eqn:reacting system 2}
  q(x,x') = \sum_{u, u' \in T, y \in S} q_{(u,y,u')}^{(\rs{D} \to \sr{A})}(x, x') + q_{\rs{i}}(x, x'), \qquad x, x' \in S,
\end{eqnarray}
where $q_{\rs{i}}$ is a transition rate function for internal transitions, and
\begin{eqnarray}
\label{eqn:q uu'}
  q_{(u,y,u')}^{(\rs{D} \to \sr{A})}(x, x') = q_{du}(x, y) r(u, y, u') p_{au'}(y, x').
\end{eqnarray}
We refer to the model described by this Markov chain as a {\bf self-reacting system}. One may think of $r$ as a routing probability matrix. For each $u, u', y$, we denote the set of pairs $(g, \sigma)$ of a transition function $g$ of the form \eqn{q uu'} for some $q_{du}$, $r$ and $p_{au'}$ and measure $\sigma$ on $S$ to support $g$ by $\sr{Q}^{(\rs{D} \to \sr{A})}_{(u,y,u')}$.

By $\nu$, we denote the solution  of the following traffic equation for each $y \in S$,
\begin{eqnarray}
\label{eqn:traffic 1}
  \sum_{u' \in T} \nu(u',y) r(u',y,u) = \nu(u,y), \qquad u \in T, y \in S.
\end{eqnarray}
We always assume that the solution $\nu$ exists and $\nu(u,y) > 0$ if $r(u',y,u) > 0$. This existence is always the case if $T$ is a finite because $\{r(u,y,u'); u, u' \in S\}$ is a transition probability matrix for each $y \in S$. However, $\nu$ may not be unique not only because any multiplication with a function of $y$ is again a solution but also because the transition probability matrix may not be irreducible. Nevertheless, it will be uniquely determined under certain reversibility in structure (see \rem{reversibility 2} below).

Let $\pi$ be a supporting measure for all the transition rate functions $q_{\rs{i}}$ and $q_{(u,y,u')}^{(\rs{D} \to \sr{A})}$ for $u,u' \in T$ and $y \in S$. Of course, the stationary distribution of $q$ is such a measure, but we do not restrict $\pi$ to be so. Define
\begin{eqnarray*}
  \tilde{q}_{(u,y,u')}^{(\rs{D} \to \sr{A})}(x, x') = \left\{\begin{array}{ll}
  \frac {\pi(x')}{\pi(x)} q_{(u',y,u)}^{(\rs{D} \to \sr{A})}(x', x), \quad & \pi(x) > 0, x, x' \in S, \vspace{1ex}\\
  0, & \mbox{otherwise},
  \end{array} \right.
\end{eqnarray*}
and
\begin{eqnarray*}
   \qquad \tilde{q}_{\rs{i}}(x, x') = \left\{\begin{array}{ll}
  \frac {\pi(x')} {\pi(x)} q_{\rs{i}}(x', x), \quad & \pi(x) > 0, x, x' \in S, \vspace{1ex}\\
  0, & \mbox{otherwise},
  \end{array} \right.
\end{eqnarray*}

Let $U = \{(u,y,u'); u, u' \in T, y \in S\}$, and let $\Gamma$ be the permutation such that $\Gamma((u,y,u')) = (u',y,u)$ for $u, u' \in T$ and $y \in S$. We like to find under what conditions we have $(\tilde{q}_{(u',y,u)}^{(\rs{D} \to \sr{A})}, \pi) \in \sr{Q}^{(\rs{D} \to \sr{A})}_{(u',y,u)}$ for all $(u,y,u') \in U$.

We first note that $(\tilde{q}_{(u',y,u)}^{(\rs{D} \to \sr{A})}, \pi) \in \sr{Q}^{(\rs{D} \to \sr{A})}_{(u',y,u)}$ is equivalent to that there are some $\tilde{q}_{du'}$, $\tilde{r}$ and $\tilde{p}_{au}$ satisfying
\begin{eqnarray}
\label{eqn:D to A 1}
  \pi(x) q_{du}(x, y) r(u, y, u') p_{au'}(y, x') = \pi(x') \tilde{q}_{du'}(x', y) \tilde{r}(u', y, u) \tilde{p}_{au}(y, x).
\end{eqnarray}

Define $\beta$, $\tilde{\beta}$, $\tilde{r}$, $\tilde{p}_{au}$ and $\tilde{q}_{au}$ as
\begin{eqnarray*}
 && \beta(u,y) = \sum_{x \in S} \pi(x) q_{du}(x,y),\qquad
  \tilde{\beta}(u',y) = \sum_{u \in T} \beta(u,y) r(u, y, u'),\\
 && \tilde{r}(u',y,u) = \frac {\beta(u,y)} {\tilde{\beta}(u',y)} r(u, y, u'),\\
 && \tilde{p}_{au}(y,x) = \frac {\pi(x)} {\beta(u,y)} q_{du}(x,y),\qquad
  \tilde{q}_{au}(x',y) = \frac {\beta(u',y)} {\pi(x')} p_{au'}(y,x').
\end{eqnarray*}
Then, we can verify \eqn{D to A 1}. That is, we always have $(\tilde{q}_{(u,y,u')}^{(\rs{D} \to \sr{A})}, \pi) \in \sr{Q}^{(\rs{D} \to \sr{A})}_{(u,y,u')}$ for all $(u,y,u') \in U$. Thus, no condition is needed for $q$ of the self-reacting system to be $\Gamma$-reversible in structure. This fact exactly explains the reversibility in (iii) of \sectn{introduction}. The problem is that the reversibility in structure gives no clue to compute $\pi$ as the stationary distribution.

To get out this situation, we consider a stronger reversibility condition. For this, we assume that the routing function $\tilde{r}$ is identical with the time reversal of $r$ with respect to its stationary measure $\nu$ of \eqn{traffic 1}. That is, we assume that
\begin{eqnarray}
\label{eqn:r reversed}
  \tilde{r}(u',y,u) = \frac {\nu(u,y)} {\nu(u',y)} r(u, y, u'), \qquad u,u' \in T, y \in S.
\end{eqnarray}
This is equivalent to $\beta = \tilde{\beta} = \nu$, which is further equivalent to
\begin{eqnarray}
\label{eqn:departure 1}
  \sum_{x' \in S} \pi(x') q_{du}(x', x) = \nu(u,x) , \qquad u \in T, x \in S.
\end{eqnarray}
We denote the subset of $\sr{Q}^{(\rs{D} \to \sr{A})}_{(u,y,u')}$ which satisfies \eqn{departure 1} for some $\nu$ of \eqn{traffic 1} by $\sr{Q}^{(\rs{D} \rightleftharpoons \sr{A})}_{(u,y,u')}$. We refer to a self-reacting system with $q_{(u,y,u')}^{(\rs{D} \to \sr{A})} \in \sr{Q}^{(\rs{D} \rightleftharpoons \sr{A})}_{(u,y,u')}$ for all $(u,y,u') \in U$ as a self-reacting system with a reversible invariant routing measure because \eqn{r reversed} is equivalent to that $\tilde{r}$ has an invariant measure $\tilde{\nu}$ such that
\begin{eqnarray*}
  \nu(u,y) = \tilde{\nu}(u,y), \qquad u \in T, x \in S.
\end{eqnarray*}

  Thus, we have the following theorem.
  
\begin{theorem} {\rm
\label{thr:reversibility 2}
  $\{\sr{Q}^{(\rs{D} \to \sr{A})}_{(u,y,u')}; (u,y,u') \in T \times S \times S\}$ is $\Gamma$-reversible in structure. Assume that the Markov chain $\{X(t)\}$ with transition rate function $q$ given by \eqn{reacting system 2} has the stationary distribution $\pi$ and $(q_{(u,y,u')}^{(\rs{D} \to \sr{A})}, \pi) \in \sr{Q}^{(\rs{D} \to \sr{A})}_{(u,y,u')}$, that is, it is a self-reacting system. Then, the time reversed Markov chain $\{\tilde{X}(t)\}$ is also a self-reacting system, and the following conditions are equivalent.
\begin{itemize}
\item [(f1)] $\{X(t)\}$ is $\Gamma$-reversible in structure with respect to $\{\sr{Q}^{(\rs{D} \rightleftharpoons \sr{A})}_{(u,y,u')}; (u,y,u') \in U\}$. That is, \eqn{departure 1} is satisfied.
\item [(f2)] $\{X(t)\}$ represents the self-reacting system with a reversible invariant routing measure. 
\item [(f3)] $\{\tilde{X}(t)\}$ is a self-reacting system with a routing function $\tilde{r}$ of \eqn{r reversed}.
\end{itemize}
}\end{theorem}
\begin{remark} {\rm
\label{rem:reversibility 2}
  $\nu$ is uniquely determined by \eqn{departure 1} if the stationary distribution $\pi$ is given.
}\end{remark}

We now construct the stationary distribution from the solution of the traffic equation \eqn{traffic 1}. From the stationary equation \eqn{stationary 2}, the stationary distribution $\pi$ must satisfy
\begin{eqnarray*}
  \pi(x) = \frac 1{a(x)} \sum_{x' \in S} \pi(x') q(x',x),
\end{eqnarray*}
where $a(x) = \sum_{x' \in S} q(x,x')$. Substituting $\nu$ of \eqn{departure 1} into this equation and using the traffic equation \eqn{traffic 1}, we have
\begin{eqnarray*}
  \pi(x) &=& \frac 1{a(x)} \sum_{x' \in S} \pi(x') \Big(\sum_{y \in S} \sum_{u, u' \in T} q_{d}(x', y, u') r(u',y,u) p_{a}(u,y, x) + q_{\rs{i}}(x', x) \Big)\\
  &=& \frac {1}{a(x)} \Big( \alpha \sum_{y \in S} \sum_{u, u' \in T} \nu(y,u') r(u',y,u) p_{a}(u,y, x) + \sum_{x' \in S} \pi(x') q_{\rs{i}}(x', x) \Big)\\
  &=& \frac {1}{a(x)} \Big( \alpha \sum_{y \in S} \sum_{u \in T} \nu(y,u) p_{a}(u,y, x) + \sum_{x' \in S} \pi(x') q_{\rs{i}}(x', x) \Big).
\end{eqnarray*}

  We assume
\begin{itemize}
\item [(e4)] $\tilde{q}_{\rs{i}}(x, x') = q_{\rs{i}}(x,x')$ for $x', x \in S$ and $b(x) \equiv \sum_{x' \in S} q_{\rs{i}}(x, x') < a(x)$ for all $x \in S$.
\end{itemize}
Then $\pi$ is uniquely determined by
\begin{eqnarray}
\label{eqn:batch movement 1}
  \pi(x) = \frac {\alpha}{a(x) - b(x)} \sum_{y \in S} \sum_{u \in T} \nu(y,u) p_{a}(u,y, x), \qquad x \in S.
\end{eqnarray}
  Thus we can construct the stationary distribution $\pi$ from $\nu$ in this case. 

\begin{example}[Balanced departure under symmetric service] {\rm
\label{exa:no internal transition 1}
  Suppose that a self-reacting system satisfies the condition (e4). We consider conditions to have \eqn{departure 1} for this network. Define $\Phi$ as
\begin{eqnarray}
\label{eqn:arrival 1}
  \Phi(x) = \frac 1{a(x)-b(x)} \sum_{y \in S} \sum_{u \in T} \nu(y,u) p_{a}(u,y, x), \qquad x \in S.
\end{eqnarray}
  We assume that state $x$ is uniquely determined by $y$ and $u$ when entity $u$ is released under state $x \in S$. Denote this $x$ by $g(y,u)$, and $q_{d}$ as
\begin{eqnarray}
\label{eqn:departure 2}
  q_{d}(x,y,u) = \frac {\nu(y,u)} {\Phi(x)} 1(x = f(y,u)), \qquad x,y \in S, u \in T.
\end{eqnarray}
  We refer to this departure rate to be {\bf $\Phi$-balanced} according to the terminology used in \cite{Serf1999} (see also \exa{State independent routing} below). Then, if $\sum_{x \in S} \Phi(x) < \infty$, we choose $\alpha = \big( \sum_{x \in S} \Phi(x) \big)^{-1}$, and let
\begin{eqnarray}
\label{eqn:batch movement 2}
  \pi(x) = \alpha \Phi(x), \qquad x \in S.
\end{eqnarray}
It is easy to check that this $\pi$ is indeed the stationary distribution and the reversibility condition \eqn{departure 1} is satisfied. Thus, by \thr{reversibility 2}, if $q_{d}$ is given by \eqn{departure 2} and if (e4) is satisfied, then we have the stationary distribution $\pi$ of \eqn{batch movement 2} as long as \eqn{arrival 1} is satisfied. This self-reacting system has a reversible invariant routing measure.
\pend
}\end{example}

\begin{example}[Bach movement network] {\rm
\label{exa:batch movement 1}
We further specialize the self-reacting system with a reversible invariant routing measure of \exa{no internal transition 1} assuming that there is no internal state transitions, that is, $b(x) \equiv 0$. Let $S = \dd{Z}_{+}^{n}$ and $T = \dd{Z}_{+}^{n+1}$, that is, its state is a $n$-dimensional vector with nonnegative integer entries. We assume that $q_{d}(x,y,u) > 0$ if and only if $y = x - u^{+}$, where $u^{+} = (u_{1}, u_{2}, \ldots, u_{n})$ for $u = (u_{0}, u_{1}, \ldots, u_{n}) \in T$. We further assume that $r(u,y,u') > 0$ if and only if $|u| = |u'|$, where $|u| = u_{0} + u_{1} + \ldots + u_{n}$.

  This system can be considered as an open queueing network with $n$ nodes and batch movements, where the $i$-th entry of $x \in S$ represents the number of customers at node $i$, and the $j$-th entry of $u \in S$ is a departing batch for $j \ne 0$ while $u_{0}$ is a departing batch from external source for $q_{d}(x,y,u) > 0$. Furthermore, this $u$ are transferred to $u'$ under each state $y$ not changing their total numbers. It is notable that, if $u_{0} = 0$ always holds, then $T$ can be reduced to $\dd{Z}^{n}$, and the model is considered as a closed network. This model is referred to as a {\bf queueing network with batch movements}.

We assume that $q_{d}(x,y,u) > 0$ only if $y = x - u^{+}$ and $p_{a}(u, y, x) > 0$ only if $x = y+u^{+}$. These assumptions are reasonable for batch movements. Then, $q_{d}(x,y,u) > 0$ if and only if $p_{a}(u, y, x) > 0$ and \eqn{departure 2} becomes
\begin{eqnarray}
\label{eqn:departure 3}
  q_{d}(x,y,u) = \frac {\nu(x-u^{+},u)} {\Phi(x)} 1(y = x-u^{+}), \qquad x,y \in S, u \in T.
\end{eqnarray}
  We show that we can choose an arbitrary $\Phi$, For this, we need to verify \eqn{arrival 1}. For the batch movement network, \eqn{arrival 1} becomes
\begin{eqnarray*}
  \Phi(x) = \frac 1{a(x)} \sum_{u \in T} \nu(x-u^{+}, u) 1(x-u^{+} \in S), \qquad x \in S.
\end{eqnarray*}
However, this automatically holds because
\begin{eqnarray*}
  a(x) &=& \sum_{x' \in S} \sum_{y \in S} \sum_{u, u' \in T} q_{d}(x, y, u) r(u,y,u') p_{a}(u',y, x') \\
  &=& \frac 1{\Phi(x)} \sum_{u \in T} \nu(x-u^{+}, u) 1(x-u^{+} \in S), \qquad x \in S.
\end{eqnarray*}
Hence, we can choose an arbitrary $\Phi$ in \eqn{departure 3} as long as it is a positive valued function. This model is called a linear network in \cite{Miya1997}. We refer to it as a {\bf batch movement network with balanced departure}. There are a lot of literature on this class of networks (see, e.g., \cite{BoucDijk1991,HendTayl1991,Serf1993}.

If all batch sizes are unit, then the batch movement network is reduced to the standard queueing network with single arrival and departures. 
\pend
}\end{example}

\begin{example}[State independent routing] {\rm
\label{exa:State independent routing}  
We consider a special case of the batch movement network with balanced departure of \exa{batch movement 1} in which $r(u,y,u')$ is independent of $y$. This model is said to have independent routing.

In this network, $\nu(y,u)$ can be independent of $y$, so it can be written as $\Psi(y) \nu(u)$ for an arbitrary positive valued function $\Psi$ on $S$ and the solution $\nu$ of the traffic equation \eqn{traffic 1} which is a function on $T$. Thus, \eqn{departure 3} becomes 
\begin{eqnarray}
\label{eqn:departure 4}
  q_{d}(x,y,u) = \frac {\Psi(x-u^{+}) \nu(u)} {\Phi(x)} 1(y = x-u^{+}), \qquad x,y \in S, u \in T,
\end{eqnarray}
  where $\Psi$ and $\Phi$ in \eqn{departure 3} can be arbitrary as long as they are positive valued functions. Furthermore, if $\nu(u)$ is of the form $\prod_{i=1}^{n} w_{i}^{u_{i}}$ for positive constants $w_{1}. w_{2}, \ldots, w_{n}$, then $\nu$ is included in $\Psi(x-u^{+}) /\Phi(x)$ since
\begin{eqnarray*}
  \nu(u) = \prod_{i=1}^{n} w_{i}^{u_{i}} = \frac {\nu(x)} {\nu(x-u^{+})}.
\end{eqnarray*}
  Thus, \eqn{departure 4} is further simplified as
\begin{eqnarray}
\label{eqn:departure 5.3}
  q_{d}(x,y,u) = \frac {\Psi(x-u^{+})} {\Phi(x)} 1(y = x-u^{+}), \qquad x,y \in S, u \in T,
\end{eqnarray}
  and, if the stationary distribution exists, then it is given by
\begin{eqnarray*}
  \pi(x) = C \Phi(x) \prod_{i=1}^{n} w_{i}^{x_{i}}, \qquad x \in S,
\end{eqnarray*}
  where $C$ is the normalizing constant. Otherwise, this $\pi$ with $C=1$ is a stationary measure. This model has been widely discussed in the literature. For example, Serfozo \cite{Serf1999} called it the {\bf Whittle network} with $\Phi$-balanced departure intensities (see Theorem 1.48 of \cite{Serf1999}). There are a lot of applications of this special case to telecommunication networks (e.g., see \cite{Bona2007}).
\pend
}\end{example}

\section{Queueing network and product form}
\label{sect:product form}

We return to reacting systems. We here use them as nodes of a queueing network with Markovian routing, where the external source is included as one of those nodes. Consider a queueing network with $n$ nodes, numbered as $1,2,\ldots,n$. We also have an external source as node $0$. Let
\begin{eqnarray*}
  J = \{ 0,1,\ldots,n \}.
\end{eqnarray*}
Each node $i$ is a reacting systems with state space $S_{i}$, departure rate $q^{(i)}_{d}$, internal transition rate $q^{(i)}_{\rs{i}}$ and arrival probability function $p^{(i)}_{a}$. Let $T_{i}$ be the set of customer types arriving at and departing from node $i$. We assume the following dynamics.
\begin{itemize}
\item [(g1)] A departure of type $u_{i} \in T_{i}$ customer from node $i$ is produced under state $x_{i} \in S_{i}$ with rate $q^{(i)}_{d u_{{i}}}(x_{i},x_{i}')$ which changes the state $x_{i}$ to $x_{i}' \in S_{i}$.
\item [(g2)] A type $u_{i}$ customer departing from node $i$ goes to node $j$ as type $u_{j}$ customer with probability $r(iu_{i},ju_{j})$ independently of everything else. It is assumed that
\begin{eqnarray*}
  \sum_{j \in J} \sum_{u_{j} \in T_{j}} r(iu_{i},ju_{j}) = 1, \qquad (i, u_{i}) \in J \times T_{i}.
\end{eqnarray*}
\item [(g3)] An arriving type $u_{j}$ customer at node $j$ under state $x_{j} \in S_{j}$ changes the sate to $x_{j}' \in S_{j}$ with probability $p^{(j)}_{au_{j}}(x_{j},x_{j}')$. It is assumed that
\begin{eqnarray*}
  \sum_{x_{j}' \in S_{j}} p^{(j)}_{au_{j}}(x_{j},x_{j}') = 1, \qquad (u_{j}, x_{j}) \in T_{j} \times S_{j}, j \in J.
\end{eqnarray*}

\item [(g4)] The state $x_{i} \in S_{i}$ at node $i$ changes to $x_{i}' \in S_{i}$ with rate $q^{(i)}_{\rs{i}}(x_{i},x_{i}')$ not producing any departure.
\end{itemize}

We model this network by a continuous time Markov chain $\{X(t)\}$. Its state space $S$ is given by
\begin{eqnarray*}
  S = S_{0} \times S_{1} \times \cdots S_{n},
\end{eqnarray*}
its transition rate function $q$ is given by
\begin{eqnarray*}
  q(x,x') = \sum_{i,j \in J} \sum_{u \in T_{i}} \sum_{v \in T_{j}} q_{iu,jv}^{(\rs{d} \to \rs{a})}(x,x') + \sum_{i \in J} q^{(i)}_{\rs{i}}(x_{i},x_{i}'),
\end{eqnarray*}
where
\begin{eqnarray}
\label{eqn:q ij}
  q_{iu,jv}^{(\rs{d} \to \rs{a})}(x,x') = \left\{\begin{array}{ll}
  q^{(i)}_{du}(x_{i}, x_{i}') r(iu,jv) p^{(j)}_{av}(x_{j},x_{j}') \prod_{k\ne i,j} 1(x_{k} = x_{k}'), \quad & i \ne j, \\
  \sum_{y_{i} \in S_{i}} q^{(i)}_{du}(x_{i}, y_{i}) r(iu,iv) p^{(i)}_{av}(y_{i},x_{i}') \prod_{k\ne i} 1(x_{k} = x_{k}'), \quad & i=j.
  \end{array} \right.
\end{eqnarray}
  We refer the model described by this transition rate as a network with reacting nodes and Markovian routing.
  
  We aim to construct a time reversed process of this network model. For this, we consider the reacting system of node $i$ with Poisson arrivals with rate $\alpha^{(i)}_{u}$. We assume:
\begin{itemize}
\item [(g5)]  There exist positive $\{\alpha^{(i)}_{u_{i}}\}$ and $\{\beta^{(i)}_{u_{i}}\}$ such that the Markov chain with transition rate $q_{i}$ given by
 \begin{eqnarray}
\label{eqn:q i}
 \hspace{-3ex} q_{i}(x_{i}, x_{i}') = \sum_{u_{i} \in T_{i}} \Big(\alpha^{(i)}_{u_{i}} p^{(i)}_{au_{i}}( x_{i},x_{i}') + q^{(i)}_{du_{i}}(x_{i}, x_{i}') \Big) + q^{(i)}_{\rs{i}}(x_{i},x_{i}') , \quad x_{i}, x_{i}' \in S_{i}, i \in J, 
\end{eqnarray}
  which has the stationary distribution $\pi_{i}$ for each $\{\alpha^{(i)}_{u_{i}}\}$, and
\begin{eqnarray}
\label{eqn:quasi-r1}
 && \sum_{x_{i}' \in S_{i}} \sum_{u_{i} \in T_{i}} \pi_{i}(x_{i}') q^{(i)}_{du_{i}}(x_{i}', x_{i}) = \beta^{(i)}_{u_{i}} \pi_{i}(x_{i}), \qquad x_{i}, x_{i}' \in S_{i}, i \in J,\\
\label{eqn:quasi-r2}
 && \alpha^{(j)}_{u_{j}} = \sum_{i \in J} \sum_{u_{i} \in T_{j}} \beta^{(i)}_{u_{i}} r(iu_{i},ju_{j}), \qquad j \in J, u_{j} \in T_{j}.
\end{eqnarray}
\end{itemize}

Then, similarly to the reacting systems with Poisson arrivals in \sectn{queue}, we define times reversed reacting system and routing function in the following way.
\begin{eqnarray*}
 && \tilde{p}_{au_{i}}^{(i)}(x_{i}, x_{i}') = \frac 1{\beta^{(i)}_{u_{i}}} \frac {\pi_{i}(x_{i}')} {\pi_{i}(x_{i})} q^{(i)}_{du_{i}}(x_{i}', x_{i}), \qquad x_{i}, x_{i}' \in S_{i}, u_{i} \in T_{i},\\
 && \tilde{q}_{du_{i}}^{(i)}(x_{i}, x_{i}') = \frac {\pi_{i}(x_{i}')} {\pi_{i}(x_{i})} \alpha^{(i)}_{u_{i}} p^{(i)}_{au_{i}}(x_{i}', x_{i}), \qquad x_{i}, x_{i}' \in S_{i}, u_{i} \in T_{i},\\
 && \tilde{q}^{(i)}_{\rs{i}}(x_{i}, x_{i}') = \frac {\pi_{i}(x_{i}')} {\pi_{i}(x_{i})} q^{(i)}_{\rs{i}}(x_{i}', x_{i}), \qquad (x_{i},x_{i}') \in S_{i},\\
 && \tilde{r}(iu_{i},ju_{j}) = \frac {\beta^{(j)}_{u_{j}}} {\alpha^{(i)}_{u_{i}}} r(ju_{j}, iu_{i}).
\end{eqnarray*}

Let $U_{i} = \{\{(a,u); u \in T_{i}\}, \{(d,u); u \in T_{i}\}\}$ for each $i \in J$, and let $U = \prod_{i \in J} U_{i}$. Define the permutation $\Gamma_{i}$ on $U_{i}$ such that $\Gamma_{i}(\{(a,u); u \in T_{i}\}) = \{(d,u); u \in T_{i}\}$ for each $i \in J$, and define the permutation on $U$ such that $\Gamma(\{v_{i}; v_{i} \in U_{i}, i \in J\}) = \{\Gamma(v_{i}); v_{i} \in U_{i}, i \in J\}$. Let $\sr{Q}^{(i)}_{au}$ and $\sr{Q}^{(i)}_{su}$ be similarly defined to $\sr{Q}_{a}$ and $\sr{Q}_{d}$, respectively, of a reacting system (see \eqn{Q a}).

The following facts can be proved using \lem{Kelly} (see \cite{ChaoMiyaPine1999} for its proof).

\begin{theorem} {\rm
\label{thr:product form 1}
  If the network with reacting nodes and Markovian routing satisfies (g5), that is, each node $i$ is quasi-reversible in separation, that is, $\Gamma_{i}$-reversible in structure with respect to $\{\sr{Q}^{(i)}_{au}, \sr{Q}^{(i)}_{du}; u \in T_{i}\}$, then it has the stationary distribution $\pi$ given by $\pi = \prod_{i=1}^{n} \pi_{i}$, that is,
\begin{eqnarray}
\label{eqn:product form 1}
  \pi(x) = \prod_{i=0}^{n} \pi_{i}(x_{i}), \qquad x \equiv (x_{0}, x_{1}, \ldots, x_{n}) \in S.
\end{eqnarray}
  Furthermore, this Markov chain is $\Gamma$-reversible in structure with respect to $\{\sr{Q}^{(\rs{D} \rightleftharpoons \sr{A})}_{(i,u), (j,v)}; u \in T_{i}, v \in T_{j}, i,j \in J\}$, where $\sr{Q}^{(\rs{D} \rightleftharpoons \sr{A})}_{(i,u), (j,v)}$ is define as the set of pair $(g,\sigma)$ of transition rate function $g$ of the form \eqn{q ij} satisfying \eqn{quasi-r1} and a measure $\sigma \equiv \prod_{i=1}^{n} \sigma_{i}$ on $S$ for $\sigma_{i}$ is a measure on $S_{i}$ such that $\sigma_{i}$ is a supporting measure of $g_{i}$ which has the same form as $q_{i}$ of \eqn{q i}.
}\end{theorem}

\begin{example} {\rm
\label{exa:product form network}
  Jackson network is well known to have the product form stationary distribution, provided it is stable. This network is a special case of the network of \thr{product form 1}. To see this, let $S = \dd{Z}_{+}^{n+1}$ and let
\begin{eqnarray*}
 && p^{(i)}_{a}(x,x') = \left\{\begin{array}{ll}
  1(x=x') & i=0,\\
  1(x_{i}' = x_{i} + 1\ge, x_{j}' = x_{j}, i \ne j), \quad & i \ne 0,
  \end{array} \right.\\
 && q^{(i)}_{d}(x,x') = \left\{\begin{array}{ll}
 \lambda & i=0,\\
  \min(x_{i},s_{i}) \mu_{i} 1(x_{i}' = x_{i} - 1\ge 0, x_{j}' = x_{j}, i \ne j), \quad & i \ne 0,
  \end{array} \right.\\
  && r(iu,jv) = p_{ij}.
\end{eqnarray*}
  Then, the reacting system with $p^{(i)}_{a}$, $q^{(i)}_{d}$ and Poisson arrivals is reversible in structure, and it is not hard to see that (g4) is satisfied with $\alpha^{(i)}_{u_{i}} = \beta^{(i)}_{u_{i}}$ for all $i$ and $u_{i}$. If we replace the $i$-th reacting systems for $i \ne 0$ in the Jackson network by a reacting systems with symmetric service, (g4) is also satisfied. This network is called a {\bf Kelly network} (see \cite{Kell1979}). Its special case limited to some symmetric service is called a {\bf BCMP network} (see \cite{BaskChanMuntPala1975}), where the workload for service to be requested is generally distributed, but any distribution can be approximated by the distribution with finitely stages.
\pend
}\end{example}

There are examples for \thr{product form 1} in other direction. A customer is said to be negative if it removes one customer without any other change at its arriving node if any, which was introduced by Gelenbe \cite{Gele1991}. Such a customer is said to be a negative signal if it instantaneously moves around in a network according to a given Markovian routing until it meets an empty node or leaves the network. It is known that \thr{product form 1} can be applied to the {\bf Jackson network with negative signals}. In this case, we may have $\alpha^{(i)}_{u_{i}} \ne \beta^{(i)}_{u_{i}}$ for some $i$ and $u_{i}$. See \cite{ChaoMiyaPine1999} for a full story of those networks.

\section{Other stochastic processes}
\label{sect:other}

We have mainly considered a continuous time Markov chain. Reversibility in structure also may be considered for other stochastic processes. A discrete time Markov chain is one of them. There are much literature for local balance, but those results are mostly parallel to those of the continuous time case. They also can be reformulated by reversibility in structure.

Both Markov chains have discrete state spaces. To directly consider service times, we need continuous valued components in the state of system. A good model for this is a generalized semi-Markov process, GSMP for short, and its variants such as reallocatable GSMP (see, e.g. \cite{Miya1993,Scha1986}). Roughly speaking, the states of those processes are the pairs of discrete state and positive real numbers, called remaining (or attaining) lifetimes. If one of the remaining lifetimes expires, then it causes the transition of the discrete state, which may produce new lifetimes.

Those GSMP's are particularly useful to prove the so-called insensitivity, which means that the stationary distribution depends on the lifetime distributions only through their means. This insensitively holds if and only if certain partial balance hold, which is a spacial case of reversibility in structure. Under this condition, we again have the symmetric queues and their networks. There are some related results here if there are service interruptions which force customers to leave (see, e.g. \cite{Miya1993}). However, no new development has been made during the last ten years.

There are some other stochastic processes related to queueing networks. For example, they are often considered using reflecting random walk or semi-Martingale reflecting Brownian motion, SRBM for short. Both are multidimensional and take values in the nonnegative orthant. Furthermore, the time reversed process has been studied for a general Markov processes (see, e.g., \cite{GetoShar1984}). It is known that SRBM on a polyhedral domain, which includes the nonnegative orthant as a special case, has the product form stationary distribution if and only if the so-called skew symmetric condition holds, which implies that the time reversed process is also SRBM (see Theorem 1.2 of \cite{Will1987a}).

\section{Concluding remarks}
\label{sect:concluding}

As we have seen, available models for queues and their networks are limited under reversibility in structure. Those models are still actively applied, but theoretical study is less active now. What can we do to gain more applicability~? There are some comments here.
\begin{enumerate}
\item The reversibility of this paper has two ingredients. One is how strongly we set up it. In this article, we mainly considered this aspect for queueing network processes. Another is which class of stochastic processes we take for considering it. We shortly discussed this issue in \sectn{other}. If we take an appropriate class, say, not too small but not too large, then we may find more tractable models for queueing networks. This direction of study seems to be not well developed.
\item The existing studies for reversibility in queues have been intended to get the stationary distribution in closed form. This target may be too big, and it may be better to have smaller target. For example, tail asymptotic of the stationary distribution is one of them (see, e.g., \cite{Miya2011}). This may include to reconsider the definition of reversibility in structure.
\end{enumerate}

\subsection*{Acknowledgements}
The author is grateful to an anonymous referee for its thoughtful comments. Reversibility in structure was greatly changed from its original form due to them.

\bibliography{miya-2013a}
 
\end{document}